\documentclass[11pt]{amsart}
\usepackage[T1]{fontenc}
\usepackage[english]{babel}
\usepackage{kpfonts}

\usepackage{amsmath,amssymb,amsthm,enumerate,xspace}
\usepackage{comment}
\usepackage{accents}
\usepackage[colorlinks=true,citecolor=blue,linkcolor=blue]{hyperref}
\usepackage{mathrsfs}
\usepackage{cleveref}

\usepackage{varwidth}

\numberwithin{equation}{section}

\newtheorem{thm}{Theorem}[section]

\newtheorem*{prob*}{Problem}

\newtheorem{prop}[thm]{Proposition}
\newtheorem{lem}[thm]{Lemma}
\newtheorem{cor}[thm]{Corollary}

\newtheorem*{iprob*}{Problem}

\theoremstyle{definition}

\newtheorem*{defi*}{Definition}

\newtheorem{exam}[thm]{Example}
\newtheorem*{acks*}{Acknowledgements}

\newcommand{\ZZ}{\mathbf{Z}}

\newcommand{\RR}{\mathbf{R}}
\newcommand{\NN}{\mathbf{N}}

\newcommand{\CC}{\mathbf{C}}
\newcommand{\QQ}{\mathbf{Q}}

\newcommand{\SO}{\mathbf{SO}}

\newcommand{\SL}{\mathbf{SL}}

\newcommand{\GL}{\mathbf{GL}}

\newcommand{\Aut}{{\rm Aut}}

\newcommand{\inv}{^{-1}}

\newcommand{\lra}{\longrightarrow}

\newcommand{\wt}{\widetilde}
\newcommand{\wh}{\widehat}

\newcommand{\ol}{\overline}

\DeclareMathOperator{\Rad}{Rad}

\newcommand{\tto}{\twoheadrightarrow}
\title[Reducing Ulam's problem to the simple case]%
{Lie groups in the symmetric group: \\
reducing Ulam's problem to the simple case}

\date{November 10, 2023}  

\author[A. Conversano]{Annalisa Conversano}
\address{Massey University Auckland, New Zealand}
\author[N. Monod]{Nicolas Monod}
\address{EPFL, Switzerland}
\begin{document}

\begin{abstract}
Ulam asked whether all Lie groups can be represented faithfully
on a countable set. We establish  a reduction of Ulam's problem to the case of simple Lie groups.

In particular, we solve the problem for all solvable Lie groups and more generally Lie groups with a linear Levi component. It follows that every amenable locally compact second countable group acts faithfully on a countable set.
\end{abstract}
\maketitle


\section{Introduction}

A group is called \textbf{countably representable} if it can be realised as a permutation group of a countable set. For obvious cardinality reasons, this only makes sense for groups that are not larger than the continuum $\mathfrak c = |\RR| = 2^{\aleph_0}$, and it is trivial for countable groups.

Therefore the rich world of \emph{Polish groups} is a prime location to study countable representability, even though (or perhaps because) this notion is non-topological. Some examples and counter-examples are discussed in~\cite{Monod_Ulam}.

Schreier and Ulam observed in 1935 that the group $\RR$ is countably representable~\cite{Scottish58}. This prompted the following problem~\cite{Ulam_book60}, still open to this day (Problem~15.8.b in~\cite{Kourovka2022}).

\begin{prob*}[Ulam]
Is every Lie group countably representable?
\end{prob*}

In 1999, Thomas~\cite[\S2]{Thomas99}, then Kallman~\cite{Kallman00} and later Ershov--Churkin~\cite{Ershov-Churkin} proved that this is the case for every \emph{linear} Lie group. There are of course many non-linear Lie groups, and although the non-linearity is in a way caused only by the center, it is well-known that the center is a fundamental obstruction to such questions, see e.g.~\cite{McKenzie71}, \cite{Churkin05} and the discussion in~\cite{Monod_Ulam}.  

Nonetheless, it was recently established that every \emph{nilpotent} Lie group is countably representable~\cite{Monod_Ulam}. Our first result applies to many non-linear groups beyond the nilpotent case: 

\begin{thm}\label{thm:intro:sol}
Every solvable Lie group is countably representable.
\end{thm}

Turning now to the most general case of Ulam's problem, recall that any Lie group $G$ admits a maximal connected solvable normal subgroup $R=\Rad(G)$, called the \textbf{radical} of $G$, which is a closed Lie subgroup. Moreover, Levi's decomposition theorem provides an \emph{immersed} connected semi-simple Lie group $S$ in $G$ such that the neutral component $G^0$ of $G$ can be written $G^0=RS$ with $R\cap S$ central in $S$, thus countable. By Malcev's theorem~\cite{Malcev42} (see~\cite[5.6]{Hilgert-Neeb}), $S$ is unique up to conjugacy and we hence abusively refer to it as ``the'' \textbf{semi-simple Levi component} of $G$. Unlike the radical, $S$ is not always closed, as witnessed by the following classical example.

 \begin{exam}
Let $a \in \SO_2(\RR)$ be a rotation through an angle incommensurable with $\pi$, and let
$b \in \wt \SL_2(\RR)$ be a generator of the center of the universal cover of $\SL_2(\RR)$. Set $H = \SO_2(\RR) \times  \wt \SL_2(\RR)$ and $\Gamma$ to be the central subgroup generated by 
$(a, b)$ in $H$. Define $G = H/\Gamma$.  

The radical $R$ of $G$ and  the semisimple Levi component $S$  are respectively the images of $\SO_2(\RR) \times \{e\}$ and $\{e\} \times \wt \SL_2(\RR)$ under the natural homomorphism $H \to H/\Gamma$. Since $\langle a \rangle$ is dense in $\SO_2(\RR)$, it follows that $S$ is not closed in $G$.
\end{exam}

The main result of this article is the following strengthening of \Cref{thm:intro:sol}:

\begin{thm}\label{thm:intro:reduction}
A Lie group is countably representable if and only if its semi-simple Levi factor is countably representable.
\end{thm}

A preliminary step in the proof is to consider the easier situation where the semi-simple Levi factor is linear, in which case it is possible to combine the proof of \Cref{thm:intro:sol} and Thomas's theorem (see \Cref{thm:levi-linear} below). This includes notably the case of compact Levi components and therefore we can apply the solution to Hilbert's fifth problem and deduce a statement far beyond Lie groups:

\begin{cor}\label{cor:intro:amen}
Every locally compact second countable group which is amenable is countably representable. 
\end{cor}

\noindent
(The corollary holds for all groups that are amenable \emph{as locally compact groups}, not only for those amenable as abstract groups.)

\medskip

In conclusion, let us consider the status of Ulam's general problem in the light of \Cref{thm:intro:reduction}. The semi-simple Levi factor can be further decomposed as a commuting product of simple Lie groups. Ulam's problem remains open because we do not know whether non-linear simple Lie groups are countably representable. Indeed we show that this is the \emph{only} obstruction:

\begin{thm}\label{thm:intro:equivalence}
The following statements are equivalent:
\begin{enumerate}[(i)]
\item every Lie group is countably representable;\label{pt:equivalence:U}
\item every connected simple Lie group with finite center is countably representable.\label{pt:equivalence:S}
\end{enumerate}
\end{thm}

\section{Preliminaries}\label{sec:prelim}

A subgroup $H<G$ of a group $G$ is said to have \textbf{countable index} (in $G$) if the coset space $G/H$ is countable. The term \emph{cocountable} is also used in the literature. We shall constantly use the basic observation that $G$ is countably representable if and only if it admits a sequence $(H_n)_{n\in \NN}$ of countable index subgroups $H_n<G$ whose intersection $\bigcap_n H_n$ is trivial.

\medskip

We recall from Lemma~11 in~\cite{Monod_Ulam} that if $G$ admits a countable index subgroup that is countably representable, then $G$ itself is countably representable.

\medskip

We shall routinely reduce ourselves from Lie groups to the case of \emph{connected} Lie groups. Indeed, since Lie groups are second countable and locally connected, they have at most countably many connected components. Therefore the connected component of the identity $G^0$ has countable index in $G$, and we can apply the above principle.

\medskip

For general locally compact groups (\Cref{cor:intro:amen}), our statement assumes second countability to ensure that the cardinality does not exceed $\mathfrak c$.

One could counter that this sufficient condition seems not to be necessary: there are locally compact groups of cardinality~$\leq\mathfrak c$ which are not second countable, e.g.\ discrete ones. However, that generality would allow for amenable groups that are not countably representable, such as McKenzie's example~\cite{McKenzie71}.

\section{The solvable case}
The goal of this section is to prove that every solvable Lie group is countably representable, as announced in \Cref{thm:intro:sol}.

The main tool for this case is as follows. Recall that polycyclic groups include notably all finitely generated nilpotent groups, see for instance~\cite[Chap.~1]{Segal83}.

\begin{prop}\label{prop:polycyclic}
Let $m$ be a positive integer and $\Gamma< \GL_m(\RR)$ a polycyclic subgroup. There exists a countable index subgroup $H< \GL_m(\RR)$ which intersects $\Gamma$ trivially.
\end{prop}

This statement will easily be reduced to the case where $\Gamma$ is cyclic, which is the object of the following lemma.

\begin{lem}\label{lem:gamma}
Let $m$ be a positive integer and $\gamma\in \SL_m(\RR)$ any element. There exists a countable index subgroup $H< \SL_m(\RR)$ which intersects trivially the cyclic subgroup generated by $\gamma$.
\end{lem}

The proof of this lemma will use a field isomorphism between $\CC$ and the field $\bigcup_{q\geq 1}\ol\QQ((t^{1/q}))$ of Puiseux series over $\ol \QQ$, as did Thomas's proof. This provides a valuation $\CC\to  \QQ \cup\{\infty\}$ with the property that the group $\SL_m(V_\CC)$ associated to the corresponding valuation ring $V_\CC$ in $\CC$ has countable index in $\SL_m(\CC)$, see Theorem~2.5 in~\cite{Thomas99}. It will be important for us to choose the field isomorphism suitably.

\begin{proof}[Proof of \Cref{lem:gamma}]
We can argue in the larger group $\SL_m(\CC)$ since all conclusions are preserved when taking the intersection of a countable index subgroup $H< \SL_m(\CC)$ with $\SL_m(\RR)$.

We choose the identification of $\CC$ with Puiseux series in the indeterminate $t$ in such a way that \itshape if $\gamma$ has any transcendental Eigenvalue, then $t$ is one such Eigenvalue. \upshape This is possible since the automorphism group of $\CC$ acts transitively on transcendentals.

Let $L_\CC<V_\CC$ be the maximal ideal and denote the corresponding congruence subgroup by $\SL_m(V_\CC; L_\CC)$. In other words, $\SL_m(V_\CC; L_\CC)$ is the kernel of the reduction morphism $\SL_m(V_\CC) \to \SL_m(\ol\QQ)$ since $\ol\QQ$ is the residue field. In particular, $\SL_m(V_\CC; L_\CC)$ also has countable index in $\SL_m(\CC)$ since $\ol\QQ$ is countable.

We now proceed to choose $H$ as a suitable conjugate of $\SL_m(V_\CC; L_\CC)$, as follows. Since $\CC$ is algebraically closed, $\gamma$ admits a Jordan normal form. Therefore, after a conjugation we may assume that $\gamma$ is in Jordan normal form and we take $H=\SL_m(V_\CC; L_\CC)$ in that conjugation.

We now claim that if $\gamma^p\in H$ for any $p\in \ZZ$, then $\gamma^p$ is the identity.

Suppose first that some Eigenvalue of $\gamma$ is transcendental. Then the same holds for $\gamma^p$ whenever $p\neq 0$. Since $\gamma$ is upper triangular, the eigenvalues of $\gamma^p$ appear as diagonal matrix coefficients. This excludes $\gamma^p\in H$ because the diagonal coefficients of elements in $H$ have valuation zero by the definition of congruence subgroups, whereas some Eigenvalue of $\gamma^p$ is $t^p$, which has valuation $p$.

We are now in the case where all Eigenvalues of $\gamma$ are algebraic. Assume $\gamma^p\in H$. It follows that all diagonal coefficients of $\gamma^p$ are~$1$. Supposing for a contradiction that $\gamma^p$ is not the identity, it follows that $\gamma$ itself was not diagonal. Consider now any non-trivial Jordan block of $\gamma$ with corresponding Eigenvalue $\lambda\in\ol\QQ$. Recall that the first subdiagonal of the Jordan block consists of~$1$s, and that consequently the corresponding coefficients on the first subdiagonal of $\gamma^p$ are $p \lambda^{p-1}$. This is a non-zero algebraic number, which contradicts $\gamma^p\in H$ because off-diagonal coefficients of elements of $H$ must be in the valuation ideal.

This contradiction proves the claim and hence the lemma.
\end{proof}

\begin{proof}[Proof of \Cref{prop:polycyclic}]
Note that in our statement $\GL_m(\RR)$ can be replaced by $\SL_m(\RR)$ since $\GL_m(\RR)$ embeds into $\SL_{m+1}(\RR)$.

By definition, $\Gamma$ admits a subnormal series
$$1 = \Gamma_0 \lhd \Gamma_1 \lhd \cdots \lhd \Gamma_d = \Gamma$$
with cyclic quotient $\Gamma_{j}/\Gamma_{j-1}$ for each $j=1, \ldots, d$. We proceed by induction on the minimal length $d$ of such a subnormal series.

The base case $d=0$ holds trivially. We now suppose that the statement has been established for all minimal lengths~$\leq d$ and consider length $d+1$. Let $\gamma$ be a generator of the cyclic group $\Gamma_1$. By \Cref{lem:gamma}, there is a countable index subgroup $H_1< \SL_m(\RR)$ intersecting $\Gamma_1$ trivially. We now consider the group $\Delta = \Gamma \cap H_1$. This is a polycyclic subgroup of $ \SL_m(\RR)$ and in fact a subnormal series witnessing polycyclicity is given by $\Delta_j = \Gamma_j \cap H_1$ since $\Delta_{j}/\Delta_{j-1}$ embeds into $\Gamma_{j}/\Gamma_{j-1}$.

However, the length $d+1$ of this series is not minimal; indeed $\Delta_1=\Delta_0$ by the choice of $H_1$. Therefore, we can apply the inductive hypothesis to $\Delta$ and obtain a countable index subgroup $H_0< \SL_m(\RR)$ intersecting $\Delta$ trivially. Now the group $H=H_0 \cap H_1$ satisfies the desired conclusion.
\end{proof}

Combining \Cref{prop:polycyclic} with Thomas's theorem for linear groups, we obtain:

\begin{thm}\label{quot:linear}
Let $m$ be a positive integer, $G< \GL_m(\RR)$ any subgroup and $\Gamma \lhd G$ a polycyclic normal subgroup of $G$.

Then the quotient group $G/\Gamma$ is countably representable.
\end{thm}

\begin{proof}
\Cref{prop:polycyclic} provides a countable index subgroup $H < \GL_m(\RR)$ which intersects $\Gamma$ trivially. Since $\GL_m(\RR)$ is countably representable (\cite[\S2]{Thomas99}), so is $H$. Now $H\cap G$ has countable index in $G$ and is countably representable. Since $H$ intersects $\Gamma$ trivially, the image $J$ of $H\cap G$ in  $G/\Gamma$ is isomorphic to $H\cap G$ and hence $J$ is countably representable. Since $J$ has countable index in $G/\Gamma$, we can apply Lemma~11 in~\cite{Monod_Ulam} and conclude that $G/\Gamma$ itself is countably representable.
\end{proof}

We can now handle all solvable Lie groups:

\begin{proof}[Proof of \Cref{thm:intro:sol}]
Let $G$ be a solvable Lie group; as explained in Section~\ref{sec:prelim}, we shall assume without loss of generality that $G$ is connected. Let $\pi  \colon \wt G \tto G$  be the universal covering map. Recall that the kernel $\ker \pi = \pi_1(G)$ is a finitely generated abelian group, see Corollary~14.2.10(iv) in~\cite{Hilgert-Neeb}. By a Theorem of Malcev, $\wt G$ is a linear group, see~\cite[2\S7]{Onishchik-Vinberg_III}. Therefore, $G$ is a quotient as those considered in \Cref{quot:linear} and it follows from that theorem that $G$ is countably representable.
\end{proof}

For future reference, we also record a straightforward bootstrap of \Cref{prop:polycyclic}:

\begin{prop}\label{prop:poly-poly}
Let $m$ be a positive integer, $G< \GL_m(\RR)$ any subgroup, $\Gamma \lhd G$ a polycyclic normal subgroup of $G$ and $\Delta < G/\Gamma$ any polycyclic subgroup of the quotient.

Then there is a countable index subgroup $H<G/\Gamma$ which meets $\Delta$ trivially.
\end{prop}

\begin{proof}
The pre-image $\wt\Delta$ of $\Delta$ in $G$ is also polycyclic. Therefore \Cref{prop:polycyclic} ensures that $G$ contains a countable index subgroup $K<G$ with $K\cap \wt\Delta$ trivial. It now suffices to take for $H$ the image of $K$ in $G/\Gamma$.
\end{proof}

\section{Linear Levi component and amenable groups}\label{sec:amen}
We recall the Levi decomposition, see~\cite[1\S4]{Onishchik-Vinberg_III} for references:

Let $G$ be any connected Lie group, let $R=\Rad(G)$ be its radical, i.e.\ the maximal connected solvable normal subgroup of $G$, which is automatically a closed Lie subgroup.

By Levi's decomposition theorem, there exists an immersed (not necessarily closed) connected semi-simple Lie group $S$ in $G$ such that $G=RS$ with $R\cap S$ countable.

\medskip

In the special case where $G$ is \emph{simply connected}, the following hold:

\begin{enumerate}[(i)]
\item both $R$ and $S$ are simply connected~\cite[\S5]{Mostow50};
\item any connected semi-simple Lie subgroup of $G$ is closed~\cite[\S6]{Mostow50}; 
\item $R\cap S$ is trivial and therefore we have a semi-direct product decomposition~\cite[1\S4]{Onishchik-Vinberg_III}.
\end{enumerate}

\begin{thm}\label{thm:levi-linear}
Let $G$ be a connected Lie group. If the semi-simple Levi component of $G$ is linear, then $G$ is countably representable. 
\end{thm}

\begin{proof}
Let $\pi \colon \wt G \tto G$ be the universal covering map of $G$ and recall that its kernel $\Gamma=\pi_1(G)$ is a finitely generated abelian group (Corollary~14.2.10(iv) in~\cite{Hilgert-Neeb}). In particular $\Gamma$ is central since $\wt G$ is connected. Denote by $\wt G = \wt R \, \wt S$ a Levi decomposition of $\wt G$. Since Levi decompositions of connected Lie groups are determined by Levi decompositions of the corresponding Lie algebra, we can assume that $\pi$ restricts to covering maps $\wt R \tto R$ and $\wt S \tto S$ for a Levi decomposition $G=R S$ of $G$. Moreover, since both $\wt R$ and $\wt S$ are simply connected as recalled above, they are in fact the universal covers of $R$ and $S$ respectively, justifying the notation.

Define $N= \Gamma \cap \wt S$. Note that $N$ is normal in $\wt G$ because $\Gamma$ is central. Hence we can consider the connected Lie group $L= \wt G/N$, which is a cover of $G$. For the reasons indicated in the first paragraph of the proof, a Levi decomposition $L=R' S'$ of $L$ is given by the images of $\wt R$ and $\wt S$ in $L$. 

We claim that $R'$ is linear. Since $\wt R \cap \wt S$ is trivial, the map $\wt R \tto R'$ is an isomorphism and hence $R'$ is simply connected. Thus the claim follows from Malcev's linearity criterion for solvable Lie groups, see~\cite[2\S7]{Onishchik-Vinberg_III}. 

Next we claim that $S'$ is linear. The kernel of the covering map $\wt S \tto S$ is the intersection between $\wt S$ and the kernel $\Gamma$ of $\wt G \tto G$. Therefore the quotient group $\wt S/N$ is isomorphic to both $S'$ 
and $S$. Since $S$ is linear by assumption, $S'$ is linear as claimed.

It now follows that $L=R' S'$ is linear by a theorem of Malcev, see~\cite[1\S5.4]{Onishchik-Vinberg_III}. Finally, since $G$ is the quotient of $L$ by the finitely generated abelian group $\Gamma / N$, we can apply \Cref{quot:linear} and deduce that $G$ is countably representable.
 \end{proof}

\begin{cor}\label{cor:Lie-amenable}
Every amenable Lie group is countably representable.
\end{cor}

\begin{proof}
Let $G$ be an amenable Lie group; we can assume $G$ connected. Let $G=RS$ be a Levi decomposition. In view of \Cref{thm:levi-linear}, it suffices to prove that $S$ is linear.

Note that the quotient $S/Z(S)$ by the center of $S$ is a quotient of $G/R$, which is amenable as a quotient of $G$. Thus $S/Z(S)$ is amenable center-free connected semi-simple Lie group, which implies that it is compact, see Theorem~1.6 in~\cite{Furstenberg63}. Since $S$ is a cover of the compact semi-simple group $S/Z(S)$, Weyl's theorem~\cite[Theorem~4.69]{Knapp02} implies that $S$ itself is compact. It is therefore linear by an application of Peter--Weyl~\cite[Corollary~4.22]{Knapp02}.
\end{proof}
 
\begin{proof}[Proof of \Cref{cor:intro:amen}]
It was exposed in Proposition~2 of~\cite{Monod_Ulam} how a positive solution to Ulam's problem would imply that \emph{every} locally compact second countable group $G$ would be countably representable. The reduction to Lie groups through the solution of Hilbert's fifth problem used in that proof (specifically: the approximation theorem as stated in \S4.6, p.~175 in~\cite{MZ55}) produces a family of Lie groups that are all quotients of an open subgroup $G_1$ of $G$. In particular, when $G$ is amenable, all those Lie groups are amenable and hence we can apply \Cref{cor:Lie-amenable} instead of relying on a hypothetical solution to Ulam's problem.
\end{proof}

\section{A reduction to simple groups}
We finally turn to the most substantial result of this article, namely that an arbitrary Lie group is countably representable if and only if its semi-simple Levi component is countably representable (\Cref{thm:intro:reduction}).

To prepare for the proof, recall that given any connected Lie group $L$, the \textbf{lineariser} of $L$, denoted by $\Lambda(L)$, is the intersection of the kernels of all finite-dimensional linear representations of $L$~\cite{Hochschild60}. The fundamental property of the lineariser (apparently due to Goto, see \cite[\S7]{Hochschild-Mostow57}) is that $L/\Lambda(L)$ is a linear Lie group, see Theorem~16.2.7 in~\cite{Hilgert-Neeb} for a proof.

\begin{lem}\label{lem:lineariser}
Let $G$ be a connected Lie group and $L$ any immersed (not necessary closed) connected Lie subgroup. Then the lineariser $\Lambda(L)$ of $L$ is central in $G$.
\end{lem}

\begin{proof}[Proof of \Cref{lem:lineariser}]
Since $G$ is connected, the adjoint representation
\[
G \lra G/Z(G) \lra \Aut(G)
\]
descends to a \emph{faithful} representation on the quotient of $G$ by its center $Z(G)$, see Lemma~9.2.21 in~\cite{Hilgert-Neeb}. However, the Lie group automorphism group $\Aut(G)$ is linear, see Theorem~1 in~\cite{Hochschild52}. Therefore, the adjoint representation is trivial on the subgroup $\Lambda(L)$, which means that $\Lambda(L)$ is contained in $Z(G)$ as claimed.
\end{proof}

We are now ready for the proof of \Cref{thm:intro:reduction}, a significant part of which is devoted to circumventing the intersection $R \cap S$ in the Levi decomposition $G=RS$.

\begin{proof}[Proof of \Cref{thm:intro:reduction}]
We can assume that our arbitrary Lie group $G$ is connected. Let $G = RS$ be a Levi decomposition of $G$. Assuming $S$ is countably representable, we will deduce that $G$ is countably representable. The other implication is obvious.

\smallskip

Set $\Lambda = \Lambda(S)$, the lineariser of $S$. By \Cref{lem:lineariser}, $\Lambda$ is central in $G$ and we can thus consider the quotient group $G/\Lambda$. Its semi-simple Levi component is the linear group $S/\Lambda$. Therefore $G/\Lambda$ is countably representable on account of \Cref{thm:levi-linear}. Hence there is a sequence of countable index subgroups $G_n < G$ whose intersection is $\Lambda$. Indeed, such a sequence can be found by taking pre-images in $G$ of a sequence of countable index subgroups of $G/\Lambda$ whose intersection is trivial.

In order to prove the theorem, it now suffices to find another sequence of countable index subgroups $H_n < G$ such that the intersection of all $H_n$ meets $\Lambda$ trivially.

\smallskip

Since $\Lambda$ is central in the semi-simple group $S$, it is a finitely generated abelian group. In particular it is countable, and therefore it suffices to find for any given $\lambda_0\neq e$ in $\Lambda$ some countable index subgroup $H < G$ avoiding $\lambda_0$. The assumption made on $S$ implies that there is a countable index subgroup $S_0 < S$ avoiding $\lambda_0$. The remainder of the proof will deal with the fact that the naive candidate $R S_0 < G$ might nonetheless contain $\lambda_0$, and therefore we need to choose a smaller subgroup $H<R S_0$.

\smallskip

Denote by $\Gamma < G$ the intersection $R \cap S$. Since $\Gamma$ is a normal solvable subgroup of the connected semi-simple group $S$, it is central in $S$. In particular, $\Gamma$ is a finitely generated abelian group.

We use the conjugation $S$-action on $R$ to form the semi-direct product $\wh G := R\rtimes S$, which is a connected Lie group with a natural quotient map $\wh G \tto G$ given by the multiplication in $G$. The kernel of this map is the group $\{(\gamma\inv, \gamma) : \gamma\in\Gamma\}$ isomorphic to $\Gamma$. Thus the pre-image $\wh \Lambda$ of $\Lambda$ in $\wh G$ is
\[
\wh \Lambda = \{(\gamma\inv, \gamma \lambda) : \gamma\in\Gamma, \lambda\in\Lambda\}.
\]
By construction this is a finitely generated metabelian group; in fact, abelian because $\Lambda$ acts trivially on $R$.

\smallskip

Viewing $\Lambda$ in $\wh G$ (as $\{e\} \times \Lambda$), we can consider the quotient group $\wh G/\Lambda = R \rtimes \left ( S/\Lambda \right )$, using again \Cref{lem:lineariser}. Since its Levi component $S/\Lambda$ is a linear group, we are again in the situation of \Cref{thm:levi-linear}. As explained in the proof of that theorem, the group $\wh G/\Lambda = R \rtimes \left ( S/\Lambda \right )$ admits a \emph{linear} (connected) cover of the form $L=\wt R \rtimes \left ( S/\Lambda \right )$. The kernel of the quotient map $L\tto \wh G/\Lambda$ is finitely generated abelian (as a subgroup of the fundamental group of $\wh G/\Lambda$, Corollary~14.2.10(iv) in~\cite{Hilgert-Neeb})) and therefore we can invoke \Cref{prop:poly-poly} to find a subgroup of countable index $J$ in $\wh G/\Lambda$ which meets trivially the image in $\wh G/\Lambda$ of the finitely generated abelian group $\wh \Lambda$.

The pre-image $\wh J$ in $\wh G$ of $J$ is a countable index subgroup of $\wh G$ such that $\wh J \cap \wh \Lambda$ is $\{e\} \times \Lambda$.

\smallskip

We claim that the image $H<G$ of the group $\wh J \cap (R\rtimes S_0) < \wh G$ under the map $\wh G \tto G$ has the desired property of avoiding $\lambda_0$.

If not, then there is $(r,s)\in \wh  J \cap (R\rtimes S_0)$ with $r s = \lambda_0$. In particular, $r = \lambda_0 s\inv$ implies $r\in \Gamma$ and then $s= r\inv  \lambda_0$ implies $(r,s) \in \wh \Lambda$. The choice of $\wh J$ now implies $(r,s) \in \{e\} \times \Lambda$. It follows $(r,s) = (e,\lambda_0)$, which contradicts $(r,s) \in R\rtimes S_0$ by the choice of $S_0$. This confirms the claim and hence completes the proof.
\end{proof}

It remains to establish \Cref{thm:intro:equivalence}. This requires two additional observations recorded in the lemmata below.

\begin{lem}\label{lem:finite-center}
Let $G$ be a connected Lie group. If every finite-sheeted cover of $G$ is countably representable, then every cover of $G$ is so too.
\end{lem}
 
\begin{proof}
Let $\wt G$ be the universal cover of $G$ and $\Gamma<Z(\wt G)$ be the fundamental group of $G$ viewed as a finitely generated abelian subgroup of $\wt G$ (Corollary~14.2.10(iv) in~\cite{Hilgert-Neeb}).

Let $\wh G$ an arbitrary cover of $G$, so that $\wh G = \wt G / \Lambda$ for some subgroup $\Lambda < \Gamma$. Since $\Gamma/\Lambda$ is finitely generated abelian, it admits a nested sequence of finite index subgroups with trivial intersection. Taking pre-images in $\Gamma$, we obtain a nested sequence of finite index subgroups $\Gamma_n < \Gamma$ with $\bigcap_n \Gamma_n = \Lambda$.

Recalling that $\Gamma$ is central since $\wt G$ is connected, we can define the groups $G_n := \wt G / \Gamma_n$, which are finite-sheeted covers of $G$ and quotients of $\wh G$. They form an inverse system with a natural morphism from $\wh G$ to the inverse limit $\varprojlim_n G_n$. This morphism is injective because $\bigcap_n \Gamma_n = \Lambda$.

Our assumption implies that each $G_n$ is countably representable. Since we realised $\wh G$ as a subgroup of an inverse limit of a \emph{countable} system of countably representable groups, it follows from Lemma~9 in~\cite{Monod_Ulam} that $\wh G$ is itself countably representable.
\end{proof}

Next, recall that every connected semi-simple group $S$ is a (not necessarily direct) product  $S=S_1 \cdots S_k$ of finitely many connected normal subgroups $S_1, \dots, S_k$ having simple Lie algebra, and such that for all $i \neq j$ the subgroups $S_i, S_j$ commute (elementwise). These subgroups $S_i$ are referred to as \textbf{the simple factors} of $S$. This decomposition follows from the decomposition of any semi-simple Lie algebra into a direct sum of simple ideals, see for instance Proposition~5.5.11 in~\cite{Hilgert-Neeb}.

\begin{lem}\label{lem:simple-factors}
 Let $S$ be a connected semi-simple Lie group and assume that the center $Z(S)$ is finite. If all the simple factors of $S$ are countably representable, then so is $S$.
\end{lem}

\begin{proof}
Let $S_1, \dots, S_k$ be the simple factors of $S$. By assumption each $S_i$ admits a sequence of countable index subgroups $S_{i,n}<S_i$ with trivial intersection.

The center $Z(S)$ is equal to the product of the centers $Z(S_i)$ since the factors commute elementwise. In particular, all simple factors have finite center. Therefore, there is no loss of generality in assuming that $S_{i,n} \cap Z(S_i)$ is trivial for all $i,n$.

The set of products $H_n = S_{1,n} \cdots S_{k,n}$ is a subgroup of $S$ since the factors $S_i$ commute elementwise, and this subgroup has countable index in $S$ for that same reason.

We claim that $H_n$ is the direct product of the groups $S_{i,n}$. Indeed, any element $z\in S_{i,n}\cap S_{j,n}$ with $i \neq j$ is centralised by $S_{j,n}$. Since the latter is dense in $S_j$, it follows that $z$ is centralised by $S_j$. But $z$ being in $S_j$, the choice of $S_{j,n}$ shows that $z$ is trivial as claimed.

It now follows that the intersection of all $H_n$ is trivial, completing the proof.
\end{proof}

\begin{proof}[End of the proof of \Cref{thm:intro:equivalence}]
Under the hypothesis that every connected simple Lie group with finite center is countably representable, we need to show that an arbitrary Lie group $G$ is countably representable. We can suppose $G$ connected and we consider a Levi decomposition $G=R S$. Thanks to \Cref{thm:intro:reduction}, it suffices to show that $S$ is countably representable.

Since $S$ is a cover of the group $S/Z(S)$, \Cref{lem:finite-center} shows that it suffices to show that every finite-sheeted cover of $S/Z(S)$ is countably representable. Such a finite-sheeted cover is a connected semi-simple Lie group with finite center, and therefore \Cref{lem:simple-factors} reduces the question to the simple factors of $S/Z(S)$. Each of these factors still has finite center and therefore is accommodated by our hypothesis.
\end{proof}

Despite our reduction of Ulam's general problem to simple factors, there is a slight disharmony between the statements of \Cref{thm:intro:reduction} and \Cref{thm:intro:equivalence}. Namely, it is unclear to the authors whether for each \emph{given} Lie group it suffices to consider the simple factors of its semi-simple Levi component.  




\bibliographystyle{amsalpha}
\bibliography{../../BIB/ma_bib}

\end{document}